\documentclass[reqno,a4paper]{amsart}
\usepackage{amssymb,url}
\usepackage[colorlinks=true,hypertex]{hyperref}
\date{This note was drafted on 19 January 2008 in Jiaozuo; Revised on 12 February 2009 in Melbourne}
\date{}
\allowdisplaybreaks[4]
\theoremstyle{plain}
\newtheorem{thm}{Theorem}

\theoremstyle{remark}
\newtheorem{rem}{Remark}
\DeclareMathOperator{\td}{d\mspace{-2mu}}

\begin{document}
\title{Logarithmic convexity of Gini means}

\author[F. Qi]{Feng Qi}
\address[F. Qi]{Research Institute of Mathematical Inequality Theory\\ Henan Polytechnic University\\ Jiaozuo City, Henan Province, 454010\\ China}
\email{\href{mailto: F. Qi <qifeng618@gmail.com>}{qifeng618@gmail.com}, \href{mailto: F. Qi <qifeng618@hotmail.com>}{qifeng618@hotmail.com}, \href{mailto: F. Qi <qifeng618@qq.com>}{qifeng618@qq.com}}
\urladdr{\url{http://qifeng618.spaces.live.com}}

\author[B.-N. Guo]{Bai-Ni Guo}
\address[B.-N. Guo]{School of Mathematics and Informatics\\ Henan Polytechnic University\\ Jiaozuo City, Henan Province, 454010\\ China}
\email{\href{mailto: B.-N. Guo <bai.ni.guo@gmail.com>}{bai.ni.guo@gmail.com}, \href{mailto: B.-N. Guo <bai.ni.guo@hotmail.com>}{bai.ni.guo@hotmail.com}}
\urladdr{\url{http://guobaini.spaces.live.com}}

\begin{abstract}
In the paper, the monotonicity and logarithmic convexity of Gini means and related functions are investigated.
\end{abstract}



\thanks{The first author was partially supported by the China Scholarship Council}

\thanks{This paper was typeset using \AmS-\LaTeX}

\maketitle

\section{Introduction}

Recall \cite{gini-mean-org} that Gini means were defined as
\begin{equation}
G(r,s;x,y)=
\begin{cases}
\biggl(\dfrac{x^s+y^s}{x^r+y^r}\biggr)^{1/(s-r)},&r\ne s;\\[1em]
\exp\biggl(\dfrac{x^r\ln x+y^r\ln y}{x^r+y^r}\biggr),&r=s;
\end{cases}
\end{equation}
where $x$ and $y$ are positive variables and $r$ and $s$ are real variables. They are also called sum mean values.
\par
There has been a lot of literature such as \cite{Czinder-pales-00, Czinder-pales-05, D-L-PMD, Farnsworth-Orr, Liu-Gini, Losonczi-JMAA, newman-pales-03, newman-sandor-03, Zs-P-G-I, sandor-gen-math-04, sandor-banach-07, 2-toader} and the related references therein about inequalities and properties of Gini means.
\par
The aim of this paper is to prove the monotonicity and logarithmic convexity of Gini means $G(r,s;x,y)$ and related functions.

\begin{thm}\label{Gini-log-conv-thm1}
Gini means $G(r,s;x,y)$ are
\begin{enumerate}
\item
increasing with respect to both $r\in(-\infty,\infty)$ and $s\in(-\infty,\infty)$;
\item
logarithmically convex with respect to both $r$ and $s$ if $(r,s)\in(-\infty,0)\times(-\infty,0)$;
\item
logarithmically concave with respect to both $r$ and $s$ if $(r,s)\in(-\infty,0)\times(-\infty,0)$.
\end{enumerate}
\end{thm}

\begin{thm}\label{Gini-log-conv-thm2}
Let
\begin{equation}
  H_{r,s;x,y}(t)=G(r+t,s+t;x,y)
\end{equation}
for $t\in\mathbb{R}$. Then Gini means $H_{r,s;x,y}(t)$ are
\begin{enumerate}
\item
increasing on $(-\infty,\infty)$;
\item
logarithmically convex on $\bigl(-\infty,-\frac{r+s}2\bigr)$;
\item
logarithmically concave on $\bigl(-\frac{r+s}2,\infty\bigr)$
\end{enumerate}
and the function
\begin{equation}\label{frak=g}
K_{r,s;x,y}(t)=H_{r,s;x,y}(t)H_{r,s;x,y}(-t)
\end{equation}
is
\begin{enumerate}
\item
increasing on $(-\infty,0)$;
\item
decreasing on $(0,\infty)$.
\end{enumerate}
\end{thm}

\begin{thm}\label{Gini-log-conv-thm3}
The function $t\mapsto t\ln H_{r,s;x,y}(t)$ is convex
\begin{enumerate}
\item
on $\bigl(-\frac{s+r}2,0\bigr)$ if $s+r>0$;
\item
on $\bigl(0,-\frac{s+r}2\bigr)$ if $s+r<0$.
\end{enumerate}
\end{thm}

\begin{rem}
The extended mean values $E(r,s;x,y)$ have properties similar to those obtained in the above theorems, see \cite{cheung-qi-mean-rgmia, cheung-qi-mean, emv-log-convex-simple.tex, Cheung-Qi-Rev.tex, exp-funct-appl-means-simp.tex}. Similar problems were also discussed in \cite{log-mean-comp-mon.tex, log-mean-comp-mon.tex-mia}.
\end{rem}

\begin{rem}
For completeness, although the monotonicity of Gini means $G(r,s;x,y)$ has been verified in \cite{Farnsworth-Orr, gini-mean-org} and related references, we would also like to give it a proof in the next section.
\end{rem}

\section{Proofs of theorems}

\begin{proof}[Proof of Theorem~\ref{Gini-log-conv-thm1}]
It is easy to see that
\begin{equation}\label{gini-mean-int-form}
\ln G(r,s;x,y)=
\begin{cases}\displaystyle
\frac1{s-r}\int_r^s\frac{x^{u}\ln x+y^{u}\ln y}{x^{u}+y^{u}}\td u,&r\ne s,\\[1em]
\dfrac{x^{r}\ln x+y^{r}\ln y}{x^{r}+y^{r}},&r=s.
\end{cases}
\end{equation}
Since
\begin{equation}\label{1-derivative}
\frac{\td}{\td u}\biggl[\frac{x^{u}\ln x+y^{u}\ln y}{x^{u}+y^{u}}\biggr]
=\frac{x^{u} y^{u} (\ln x-\ln y)^2} {(x^{u}+y^{u})^2}>0
\end{equation}
and
\begin{equation}\label{2-derivative}
\frac{\td{}^2}{\td u^2}\biggl[\frac{x^{u}\ln x+y^{u}\ln y}{x^{u}+y^{u}}\biggr]
=-\frac{x^uy^u(x^u-y^u) (\ln x-\ln y)^3}{(x^u+y^u)^3}
=\begin{cases}
\ge0,&u\le0,\\
\le0,&u\ge0,
\end{cases}
\end{equation}
then the integrand in \eqref{gini-mean-int-form} is increasing on $u\in(-\infty,\infty)$, logarithmically convex on $u\in(-\infty,0)$ and logarithmically concave on $u\in(0,\infty)$. It is known \cite[Lemma~1]{qx3} that if $f(t)$ is differentiable and increasing on an interval $I$, then the integral arithmetic mean of $f(t)$,
\begin{equation}
\phi(r,s)=
\begin{cases}\displaystyle
\frac{1}{s-r}\int_r^sf(t)\td t,& r\ne s,\\
f(r),& r=s,
\end{cases}
\end{equation}
is also increasing with $r$ and $s$ on $I$; If $f(t)$ is twice differentiable and convex on $I$, then $\phi(r,s)$ is also convex with $r$ and $s$ on $I$. Consequently, Gini means $G(r,s;x,y)$ with respect to both $r$ and $s$ are increasing on $(-\infty,\infty)$, logarithmically convex if $(r,s)\in(-\infty,0)\times(-\infty,0)$, and logarithmically concave if $(r,s)\in(0,\infty)\times(0,\infty)$. The proof of Theorem~\ref{Gini-log-conv-thm1} is complete.
\end{proof}

\begin{proof}[Proof of Theorem~\ref{Gini-log-conv-thm2}]
Taking the logarithm of $H_{r,s;x,y}(t)$ and differentiating consecutively yields
\begin{align*}
\ln H_{r,s;x,y}(t)&=\frac1{s-r}\bigl[\ln\bigl(x^{s+t}+y^{s+t}\bigr) -\ln\bigl(x^{r+t}+y^{r+t}\bigr)\bigr],\\
\bigl[\ln H_{r,s;x,y}(t)\bigr]'&=\frac1{s-r}\biggl(\frac{x^{s+t}\ln x+y^{s+t}\ln y}{x^{s+t}+y^{s+t}} -\frac{x^{r+t}\ln x+y^{r+t}\ln y}{x^{r+t}+y^{r+t}}\biggr),\\
\bigl[\ln H_{r,s;x,y}(t)\bigr]''&=\frac1{s-r}\biggl[\frac{x^{s+t} y^{s+t} (\ln x-\ln y)^2} {(x^{s+t}+y^{s+t})^2} -\frac{x^{r+t} y^{r+t} (\ln x-\ln y)^2}{(x^{r+t}+y^{r+t})^2}\biggr].
\end{align*}
\par
By virtue of \eqref{1-derivative}, it easily follows that $\bigl[\ln H_{r,s;x,y}(t)\bigr]'\ge0$, which means that Gini means $H_{r,s;x,y}(t)$ is increasing on $(-\infty,\infty)$.
\par
With the aid of \eqref{2-derivative}, it may be obtained easily that the function
\begin{equation}
f_{x,y}(u)=\frac{x^{u} y^{u} (\ln x-\ln y)^2} {(x^{u}+y^{u})^2}
\end{equation}
is increasing on $u\in(-\infty,0)$ and decreasing on $(0,\infty)$. Moreover, it is clear that $f_{x,y}(u)=f_{x,y}(-u)$, that is, the function $f_{x,y}(u)$ is even on $(-\infty,\infty)$.
\par
Let
\begin{equation}
  F_{x,y}(t)=f_{x,y}(s+t)-f_{x,y}(r+t).
\end{equation}
If $s+t>r+t>0$, that is, $t>-r>-s$, since $f_{x,y}(u)$ is decreasing on $(0,\infty)$, then $F_{x,y}(t)\le0$. Similarly, if $r+t<s+t<0$, i.e., $t<-s<-r$, then $F_{x,y}(t)\ge0$. If $r+t<0<s+t$ and $0<-(r+t)<s+t$, equivalently, $t>-\frac{r+s}2$, since $f_{x,y}(u)$ is even on $(-\infty,\infty)$ and decreasing on $(0,\infty)$, then $F_{x,y}(t)\le0$; Similarly, if $t<-\frac{r+s}2$, then $F_{x,y}(t)\ge0$. This implies
\begin{equation*}
  \bigl[\ln H_{r,s;x,y}(t)\bigr]''
  \begin{cases}
    \ge0,&t<-\dfrac{r+s}2\\[0.5em]
    \le0,&t>-\dfrac{r+s}2
  \end{cases}
\end{equation*}
for all $r,s,x,y$ by a recourse to symmetric properties $G(r,s;x,y)=G(s,r;x,y)=G(r,s;y,x)$.
\par
\par
Taking the logarithm on both sides of \eqref{frak=g} and differentiating gives
\begin{equation*}
\bigl[\ln K_{r,s;x,y}(t)\bigr]' =\frac{H_{r,s;x,y}'(t)}{H_{r,s;x,y}(t)} -\frac{H_{r,s;x,y}'(-t)}{H_{r,s;x,y}(-t)},\\
\end{equation*}
where
\begin{equation*}
H_{r,s;x,y}'(-t)=\frac{\td H_{r,s;x,y}(u)}{\td u}\biggl|_{u=-t}.
\end{equation*}
The logarithmic convexities of $H_{r,s;x,y}(t)$ implies that the function $\frac{H_{r,s;x,y}'}{H_{r,s;x,y}}(t)$ is increasing on $\bigl(-\infty,-\frac{r+s}2\bigr)$ and decreasing on $\bigl(-\frac{r+s}2,0\bigr)$. Careful computation can verify that
\begin{align*}
\frac{H_{r,s;x,y}'(t)}{H_{r,s;x,y}(t)} =\frac{H_{r,s;x,y}'(-t-(s+r))}{H_{r,s;x,y}(-t-(s+r))}
\end{align*}
for $t\in(-\infty,\infty)$. Consequently, the function
\begin{equation*}
Q(t)=\frac{H_{r,s;x,y}'(t-(s+r)/2)}{H_{r,s;x,y}(t-(s+r)/2)}
\end{equation*}
is increasing on $(-\infty,0)$ and decreasing on $(0,\infty)$ and satisfies $Q(t)=Q(-t)$ for $t\in(-\infty,\infty)$. Utilization of the approach applied to the function $f_{x,y}$ above yields that $Q(t+(s+r))-Q(t)$ is positive on $\bigl(-\infty,-\frac{s+r}2\bigr)$ and negative on $\bigl(-\frac{s+r}2,\infty\bigr)$, which is equivalent to
\begin{equation}\label{q-q-w}
Q\biggl(t+\frac{s+r}2\biggr)-Q\biggl(t-\frac{s+r}2\biggr) =\frac{H_{r,s;x,y}'(t)}{H_{r,s;x,y}(t)} -\frac{H_{r,s;x,y}'(t-(s+r))}{H_{r,s;x,y}(t-(s+r))}
\end{equation}
being positive on $(-\infty,0)$ and negative on $(0,\infty)$. Since
\begin{equation}
K_{r,s;x,y}(t)=\frac{xyH_{r,s;x,y}(t)}{H_{r,s;x,y}(t-(s+r))},
\end{equation}
then the function in \eqref{q-q-w} equals $\bigl[\ln K_{r,s;x,y}(t)\bigr]'$, which implies that the function $K_{r,s;x,y}(t)$ is increasing on $(-\infty,0)$ and decreasing on $(0,\infty)$. The proof of Theorem~\ref{Gini-log-conv-thm2} is complete.
\end{proof}

\begin{proof}[Proof of Theorem~\ref{Gini-log-conv-thm3}]
Direct calculation yields
\begin{equation}\label{klm}
\bigl[t\ln H_{r,s;x,y}(t)\bigr]''=2\bigl[\ln H_{r,s;x,y}(t)\bigr]' +t\bigl[\ln H_{r,s;x,y}(t)\bigr]''.
\end{equation}
By Theorem~\ref{Gini-log-conv-thm2}, it follows that $\bigl[\ln H_{r,s;x,y}(t)\bigr]'>0$ on $(-\infty,\infty)$, $\bigl[\ln H_{r,s;x,y}(t)\bigr]''>0$ on $\bigl(-\infty,-\frac{s+r}2\bigr)$ and $\bigl[\ln H_{r,s;x,y}(t)\bigr]''<0$ on $\bigl(-\frac{s+r}2,\infty\bigr)$. Therefore, if $s+r<0$ then $\bigl[t\ln H_{r,s;x,y}(t)\bigr]''>0$ and $t\ln H_{r,s;x,y}(t)$ is convex on $\bigl(0,-\frac{s+r}2\bigr)$, if $s+r>0$ then $\bigl[t\ln H_{r,s;x,y}(t)\bigr]''>0$ and $t\ln H_{r,s;x,y}(t)$ is convex on $\bigl(-\frac{s+r}2,0\bigr)$. The proof of Theorem~\ref{Gini-log-conv-thm3} is complete.
\end{proof}


\begin{thebibliography}{99}

\bibitem{cheung-qi-mean-rgmia}
W.-S. Cheung and F. Qi, \textit{Logarithmic convexity of the one-parameter mean values}, RGMIA Res. Rep. Coll. \textbf{7} (2004), no.~2, Art.~15, 331\nobreakdash--342; Available online at \url{http://www.staff.vu.edu.au/rgmia/v7n2.asp}.

\bibitem{cheung-qi-mean}
W.-S. Cheung and F. Qi, \textit{Logarithmic convexity of the one-parameter mean values}, Taiwanese J. Math. \textbf{11} (2007), no.~1, 231\nobreakdash--237.

\bibitem{Czinder-pales-00}
P. Czinder and Zs. P\'ales, \textit{A general Minkowski-type inequality for two variable Gini means}, Publ. Math. Debrecen \textbf{57} (2000), no.~1-2, 203\nobreakdash--216.

\bibitem{Czinder-pales-05}
P. Czinder and Zs. P\'ales, \textit{Local monotonicity properties of two-variable Gini means and the comparison theorem revisited}, J. Math. Anal. Appl. \textbf{301} (2005), no.~2, 427\nobreakdash--438.

\bibitem{D-L-PMD}
Z. Dar\'oczy and L. Losonczi, \emph{\"Uber den Vergleich von Mittelwerten}, Publ. Math. Debrecen \textbf{17} (1970), 289\nobreakdash--297.

\bibitem{Farnsworth-Orr}
D. Farnsworth and R. Orr, \textit{Gini means}, Amer. Math. Monthly \textbf{93} (1986), no.~8, 603\nobreakdash--607.

\bibitem{gini-mean-org}
C. Gini, \textit{Di una formula compresive delle medie}, Metron \textbf{13} (1938), 3\nobreakdash--22.

\bibitem{emv-log-convex-simple.tex}
B.-N. Guo and F. Qi, \textit{A simple proof of logarithmic convexity of extended mean values}, Numer. Algorithms (2009), in press; Available online at \url{http://dx.doi.org/10.1007/s11075-008-9259-7}.

\bibitem{Liu-Gini}
Zh. Liu, \textit{Minkowski's inequality for extended mean values}, Proceedings of the Second ISAAC Congress, Vol. 1 (Fukuoka, 1999), 585--592, Int. Soc. Anal. Appl. Comput., 7, Kluwer Acad. Publ., Dordrecht, 2000.

\bibitem{Losonczi-JMAA}
L. Losonczi, \emph{Inequalities for integral mean values}, J. Math. Anal. Appl. \textbf{61} (1977), 586--606.

\bibitem{newman-pales-03}
E. Neuman and Zs. P\'ales, \textit{On comparison of Stolarsky and Gini means}, J. Math. Anal. Appl. \textbf{278} (2003), no.~2, 274\nobreakdash--284.

\bibitem{newman-sandor-03}
E. Neuman and J. S\'andor, \textit{Inequalities involving Stolarsky and Gini means}, Math. Pannon. \textbf{14} (2003), no.~1, 29\nobreakdash--44.

\bibitem{Zs-P-G-I}
Zs. P\'ales, \emph{Comparison of two variable homogeneous means, General Inequalities 6. Proc. 6th
Internat. Conf. Math. Res. Inst. Oberwolfach}, Birkh\"auser Verlag Basel, 1992, pp.~59--69.

\bibitem{log-mean-comp-mon.tex}
F. Qi, \textit{Complete monotonicity of logarithmic mean}, RGMIA Res. Rep. Coll. \textbf{10} (2007), no.~1, Art.~18; Available online at \url{http://www.staff.vu.edu.au/rgmia/v10n1.asp}.

\bibitem{log-mean-comp-mon.tex-mia}
F. Qi and Sh.-X. Chen, \textit{Complete monotonicity of the logarithmic mean}, Math. Inequal. Appl. \textbf{10} (2007), no.~4, 799\nobreakdash--804.

\bibitem{Cheung-Qi-Rev.tex}
F. Qi, P. Cerone, S. S. Dragomir and H. M. Srivastava, \textit{Alternative proofs for monotonic and logarithmically convex properties of one-parameter mean values}, Appl. Math. Comput. \textbf{208} (2009), no.~1, 129\nobreakdash--133; Available online at \url{http://dx.doi.org/10.1016/j.amc.2008.11.023}.

\bibitem{exp-funct-appl-means-simp.tex}
F. Qi and B.-N. Guo, \textit{The function $(b^x-a^x)/x$: Logarithmic convexity and applications to extended mean values}, Available online at \url{http://arxiv.org/abs/0903.1203}.

\bibitem{qx3}
F. Qi, S.-L. Xu and L. Debnath, \textit{A new proof of monotonicity for extended mean values}, Internat. J. Math. Math. Sci. \textbf{22} (1999), no.~2, 415\nobreakdash--420.

\bibitem{sandor-gen-math-04}
J. S\'andor, \textit{A note on the Gini means}, Gen. Math. \textbf{12} (2004), no.~4, 17\nobreakdash--21.

\bibitem{sandor-banach-07}
J. S\'andor, \textit{The Schur-convexity of Stolarsky and Gini means}, Banach J. Math. Anal. \textbf{1} (2007), no.~2, 212\nobreakdash--215.

\bibitem{2-toader}
S. Toader and G. Toader, \textit{Complementaries of Greek means with respect to Gini means},  Int. J. Appl. Math. Stat. \textbf{11} (2007), no.~7, 187\nobreakdash--192.

\end{thebibliography}
\end{document}